\documentclass[12pt]{amsart}
\usepackage{bm}

\newcommand{\newsection}[1]{\setcounter{equation}{0} \section{#1}}
\newcommand{\bea}{\begin{eqnarray}}
\newcommand{\eea}{\end{eqnarray}}

\newcommand{\cld}{\mathcal{D}}
\newcommand{\cle}{\mathcal{E}}

\newcommand{\clh}{\mathcal{H}}
\newcommand{\clk}{\mathcal{K}}
\newcommand{\cll}{\mathcal{L}}
\newcommand{\clm}{\mathcal{M}}

\newcommand{\clo}{\mathcal{O}}

\newcommand{\cls}{\mathcal{S}}
\newcommand{\clt}{\mathcal{T}}

\newcommand{\clz}{\mathcal{Z}}

\newcommand{\raro}{\rightarrow}

\def \qed {\hfill \vrule height6pt width 6pt depth 0pt}
\def\textmatrix#1&#2\\#3&#4\\{\bigl({#1 \atop #3}\ {#2 \atop #4}\bigr)}
\def\dispmatrix#1&#2\\#3&#4\\{\left({#1 \atop #3}\ {#2 \atop #4}\right)}
\newcommand{\be}{\begin{equation}}
\newcommand{\ee}{\end{equation}}
\newcommand{\ben}{\begin{eqnarray*}}
\newcommand{\een}{\end{eqnarray*}}

\newcommand{\NI}{\noindent}

\newcommand{\bi}{\begin{itemize}}
\newcommand{\ei}{\end{itemize}}

\def\5{{5\superprime}}

\newtheorem{Theorem}{\sc Theorem}[section]
\newtheorem{Lemma}[Theorem]{\sc Lemma}
\newtheorem{Proposition}[Theorem]{\sc Proposition}
\newtheorem{Corollary}[Theorem]{\sc Corollary}
\newtheorem{Definition}[Theorem]{\sc Definition}
\newtheorem{Example}[Theorem]{\sc Example}
\newtheorem{Remark}[Theorem]{\sc Remark}
\newtheorem{Note}[Theorem]{\sc Note}
\newtheorem{Question}{\sc Question}
\newtheorem{ass}[Theorem]{\sc Assumption}
\newcommand{\bt}{\begin{Theorem}}
\def\beginlem{\begin{Lemma}}
\def\beginprop{\begin{Proposition}}
\def\begincor{\begin{Corollary}}
\def\begindef{\begin{Definition}}
\def\beginexamp{\begin{Example}}
\def\beginrem{\begin{Remark}}
\def\beginq{\begin{Question}}
\def\beginass{\begin{ass}}
\def\beginnote{\begin{Note}}
\newcommand{\et}{\end{Theorem}}
\def\endlem{\end{Lemma}}
\def\endprop{\end{Proposition}}
\def\endcor{\end{Corollary}}
\def\enddef{\end{Definition}}
\def\endexamp{\end{Example}}
\def\endrem{\end{Remark}}
\def\endq{\end{Question}}
\def\endass{\end{ass}}
\def\endnote{\end{Note}}

\begin{document}

\title{Essentially Reductive Weighted Shift Hilbert Modules}

\author[Douglas]{Ronald G. Douglas}
\author[Sarkar]{Jaydeb Sarkar}

\address{Texas A \& M University, College Station, Texas 77843, USA}

\email{rdouglas@math.tamu.edu, jsarkar@math.tamu.edu}

\thanks{This research was partially supported by a grant from the National Science Foundation.}

\keywords{Hilbert modules, spherical isometries, weighted shifts}

\subjclass[2000]{46E22, 46M20, 47A13, 47B32}

\begin{abstract}
We discuss the relation between questions regarding the essential normality of finitely generated essentially spherical isometries and some results and conjectures of Arveson and Guo-Wang on the closure of homogeneous ideals in the $m$-shift space. We establish a general results for the case of two tuples and ideals with one dimensional zero variety. Further, we show how to reduce the analogous question for quasi-homogeneous ideals, to those results for homogeneous ones. Finally, we show that the essential reductivity of positive regular Hilbert modules is directly related to a generalization of the Arveson problem.
\end{abstract}

\maketitle

\newsection{Introduction}

\vspace{1.5cm}

Not all isometries on a complex Hilbert space are unitary or even essentially unitary; that is, unitary modulo the compacts. A unilateral shift of infinite multiplicity is a counter-example. However, if the isometry $V$ has a finite generating set, then $V$ is essentially unitary or, equivalently in this case, essentially normal.

What if the operator is only essentially isometric or $T^*T - I$ is compact and $T$ has a finite generating set? The answer is still affirmative.

\vspace{0.2in}

\begin{Theorem}\label{berger} 
If $T$ is an essentially isometric operator with a finite generating set, then $T$  is essentially unitary.
\end{Theorem}

\vspace{0.1in}

The assumption that $I - T^* T$ is compact implies that the range of $T^*$ is closed and has finite co-dimension. Thus $T$ is left semi-Fredholm. Moreover, the fact that $T$ has a finite generating set implies that the range of $T$ has finite co-dimension. Therefore $T$ is Fredholm which yields the result. 

We want to consider the possible validity of analogues of this result for commuting $m$-tuples of operators on a complex Hilbert space $\clh$ and their relation to same conjectures and results of Arveson, Guo, Wang and the first author.

An $m$-tuple of operators $(T_1, T_2, \ldots , T_m)$ on $\clh$ is said to be a spherical isometry if $$\sum_{i = 1}^{m} T_i^* T_i = I_{\clh} \;\;\;\;\;\;\;\;\;\; \mbox{or} \;\;\;\;\;\;\;\; \|f\|^2= \sum_{i=1}^{m} \|T_i f\|^2$$
for $f$ in $\clh$.

Examples of such $m$-tuples of operators are those defined to be the coordinate multiplication operators on the Bergman and Hardy spaces over the unit ball $\mathbb{B}^m$ in $\mathbb{C}^m$. On the $m$-shift space, $H^2_m$, the corresponding $m$-tuple is an essentially spherical isometry in the sense that the operator $I_{\clh} - \sum_{i = 1}^{m} T_i^* T_i$ is compact. If one takes an infinite direct sum of either of the first two examples of $m$-tuples, one obtains a spherical isometry that is not an essentially spherical unitary or $I_{\clh} - \sum_{i = 1}^{m} T_i T_i^*$ is not compact.

Although the definition of spherical isometry does not require the operators $\{T_i\}$ to commute, we will make that assumption from now on and consider various questions related to Theorem \ref{berger} in this context. In particular, if $(T_1, T_2, \ldots, T_m)$ is a commuting essentially spherical isometry on the complex Hilbert space $\clh$ which has a finite generating set, must $I_{\clh} - \sum_{i=1}^{m} T_i T^*_i$ be compact?

This question has an easy negative answer which one can see by setting $T_i = \frac{1}{\sqrt{m}} M_{z_i}$ for $i = 1, 2, \ldots, m$, where $M_{z_i}$ is coordinate multiplication by $z_i$ on the Hardy space $H^2(\mathbb{D}^m)$ on the polydisk $\mathbb{D}^m$. 

One source of the difficulty in proceeding from one variable to several can be seen by considering the notion of left semi-Fredholmness for commuting $m$-tuples. In particular, the $m$-tuple $(T_1, T_2, \ldots, T_m)$ on $\clh$ is said to be left semi-Fredholm if $\mbox{dim} \clh/ \{ T_1^* \clh + T_2^* \clh + \cdots + T^*_m \clh\} < \infty$. This implies that $T_1^* \clh + T_2^* \clh + \cdots + T^*_m \clh$ is closed and has finite co-dimension in $\clh$. 

In the one-variable case, $I - T^* T$ compact implies not only that $T$ is left semi-Fredholm but that the same is true for $T- \lambda$ for $\lambda$ in the open unit disk $\mathbb{D}$. The analogous statement fails in the case of several variables; that is, while $I_{\clh} - \sum_{i=1}^{m} T^*_i T_i$ compact implies $(T_1, T_2, \ldots, T_m)$ is left semi-Fredholm, that is not necessarily the case for $( T_1 - \lambda_1, T_2 - \lambda_2, \ldots,  T_m - \lambda_m)$ for $\bm{\lambda} = (\lambda_1, \lambda_2, \ldots, \lambda_m)$ in $\mathbb{B}^m$. The preceding example on the polydisk provides an example of this behavior. Hence we add that assumption to obtain:

\vspace{0.2in}

\begin{Question}\label{Q1} If $(T_1, T_2, \ldots, T_m)$ is a commuting essentially spherical isometry on the complex Hilbert space $\clh$ which has a finite generating set and for which $( T_1 - \lambda_1, T_2 - \lambda_2, \ldots,  T_m - \lambda_m)$ is left semi-Fredholm for $(\lambda_1, \lambda_2, \ldots, \lambda_m)$ in $\mathbb{B}^m$, then must $(T_1, T_2, \ldots, T_m)$ be an essential unitary? What if one assumes, in addition, that the operators $\{T_i\}$ are hyponormal? Or, jointly hyponormal?
\end{Question}

\vspace{0.1in}

These question are related to studies of Eschmeier and Putinar \cite{EP} and Gleason, Richter and Sundberg \cite{GRS}. In the first note, the authors survey some results on spherical isometries and present an interesting example which we will discuss in Section 7. In the latter paper, the authors discuss some examples which demonstrate the necessity of the assumption that the $m$-tuple has a finite generating set.

 Note that if the $m$-tuple is actually a spherical isometry, then both additional assumptions on hyponormality follow by Athavale's result \cite{A}. Working modulo the compacts we see that an essentially spherical isometry is essentially subnormal and hence essentially jointly hyponormal. 

Note that in the context of Question \ref{Q1}, $I_{\clh} - \sum_{i=1}^{m} T_i T_i^*$ is compact if and only if each of the $T_i$ are essentially normal.

One can rephrase these questions in the language of Hilbert modules over the algebra $\mathbb{C}[{\mathbf{z}}]$ of polynomials in $m$ variables with $\bm{z} = (z_1, z_2, \ldots, z_m)$. We will use $(M_{z_1}, M_{z_2}, \ldots, M_{z_m})$  to denote the $m$-tuple of operators defined on a Hilbert module $\clh$ by module multiplication by $z_1, z_2, \ldots, z_m$, respectively. The Hilbert module $\clh$ is said to be {\it isometric} ({\it essentially isometric} or {\it essentially unitary}) if the $m$-tuple  $(M_{z_1}, M_{z_2}, \ldots, M_{z_m})$ is a commuting spherical isometry (essentially spherical isometry or essentially spherical unitary). Then the above questions can be rephrased as follows: 

\vspace{0.2in}

\begin{Question}\label{Q2} Is every finitely generated essentially isometric Hilbert module over $\mathbb{C}[\mathbf{z}]$ for which $(M_{z_1} - \lambda_1, M_{z_2} - \lambda_2, \ldots, M_{z_m} - \lambda_m)$ is left semi-Fredholm for $(\lambda_1, \lambda_2, \ldots, \lambda_m)$ in $\mathbb{B}^m$ necessarily essentially unitary? What if the operators $\{M_{z_1}, M_{z_2}, \ldots, M_{z_m}\}$ are hyponormal?
\end{Question}

\vspace{0.2in}

We will show in this note that an affirmative answer to Question \ref{Q1} or \ref{Q2} implies an affirmative answer to Arveson's conjecture (\cite{A2}, \cite{A1}, \cite{A5}, \cite{A7}) concerning the closure of homogeneous ideals in the $m$-shift space and the Guo-Wang conjecture (\cite{GW2}, \cite{GW}) concerning ideals which are quasi-homogeneous. 

There are a couple of interesting questions in the cyclic case one can formulate by making additional assumptions. 

The Hilbert module $\clh$ over $\mathbb{C}[\mathbf{z}]$ is said to be a {\it weighted shift} Hilbert module if there is a wandering cyclic vector relative to the monomials $\{\mathbf{z}^{\bm{\alpha}}\}$; that is, there is a vector $f$ in $\clh$ such that $\{\mathbf{z}^{\bm{\alpha}} f\}$ is an orthogonal basis for $\clh$ with $\bm{\alpha} = (\alpha_1, \alpha_2, \ldots, \alpha_m)$ in $\mathbb{N}^m$ and $\bm{z}^{\bm{\alpha}} = z_1^{\alpha_1} z_2^{\alpha_2} \cdots z_m^{\alpha_m}$. The vector $f$ in $\clh$ will be said to be {\it weakly wandering} if $\{\mathbf{z}^{\bm{\alpha}} f\}$ is orthogonal to $\mathbf{z}^{\bm{\beta}}f$ if $\alpha_1 + \alpha_2 + \cdots + \alpha_m = |\alpha| \neq |\beta| = \beta_1 + \beta_2 + \cdots \beta_m$, or their degrees are distinct.

\vspace{0.2in}

\begin{Question}\label{Q3} If a $\mathbb{C}[\mathbf{z}]$ Hilbert module action on $\clh$ is essentially isometric and has a wandering (or weakly wandering) cyclic vector and $(M_{z_1} - \lambda_1, M_{z_2} - \lambda_2, \ldots, M_{z_m} - \lambda_m)$ is left semi-Fredholm for $(\lambda_1, \lambda_2, \ldots, \lambda_m)$ in $\mathbb{B}^m$, must $\clh$ be essentially unitary?
\end{Question}

\vspace{0.2in}

An affirmative answer to the latter question assuming the existence of a weakly wandering cyclic vector would yield affirmative answers to the Arveson and Guo-Wang conjectures for the closure of homogeneous and quasi-homogeneous ideals in $H^2_m$. Moreover, the $m$ in this question corresponds to the $m$ in the conjectures.  

In Sections 2 and 3, we study these questions and their relation to the various conjectures including an introduction of a family of Hilbert modules over $\mathbb{C}[\bm{z}]$ sharing many properties of the Bergman, Hardy and $m$-shift spaces.

 In Section 4 we establish affirmative answers to question \ref{Q3} in the $m=2$ case for this family of Hilbert modules and, hence, extend solutions to the conjectures of Arveson and Guo-Wang beyond that for ideals in $H^2_2$ proved by Guo and Wang (\cite{GW}, \cite{GW2}). 

In Section 5, we extend the results of Guo and Wang on submodules defined as the closure of homogeneous ideals with a one dimensional zero variety. 

In Section 6 we show that the essential reductivity of Hilbert modules over $\mathbb{C}[\rm{z}]$ defined by a positive regular polynomial in $m$ variables is equivalent to that for certain related quasi-homogeneous ideals in the $(m+k)$-shift space, where $m+k$ is the numbers of monomials in the polynomial with non-zero coefficients. We conclude in Section 7 with some additional comments on these questions and their possible resolution.

\vspace{.6 in}

\newsection{The Basic Setup}

\vspace{.4in}

Let $\clh$ be a Hilbert space completion of the polynomials $\mathbb{C}[\mathbf{z}]$, where $\mathbf{z} = (z_1, z_2, \ldots, z_m)$ for a positive integer $m$, such that each operator $M_p$ is bounded on $\clh$, where $M_p$ is defined to be multiplication by the polynomial $p(\mathbf{z})$ in $\mathbb{C}[\mathbf{z}]$. These assumptions make $\clh$ into a Hilbert module over $\mathbb{C}[\mathbf{z}]$ (cf. \cite{D-P}).

Standard examples of such Hilbert modules are given by the Hardy and Bergman spaces for the unit ball $\mathbb{B}^m$ in $\mathbb{C}^m$, the $m$-shift or symmetric Fock space in $m$ variables, or the Bergman space for certain Reinhardt domains in $\mathbb{C}^m$. Another class of examples, based on $m$- commuting weighted shifts, is discussed in \cite{D}. These examples all have property

\vspace{0.1in}

{\bf (A)}  the monomials  $\{\mathbf{z}^{\bm{\alpha}}\}$ for  $\bm{\alpha} = (\alpha_1, \alpha_2, \ldots \alpha_m)$ in $\mathbb{N}^m$ are orthogonal.

\vspace{0.1in}

We refer to (see \cite{D}) a Hilbert space completion $\clh$ of $\mathbb{C}[\mathbf{z}]$ satisfying {\bf (A)} as a {\it weighted shift} Hilbert module. All the examples mentioned above satisfy {\bf (A)}. 

A monomial $\mathbf{z}^{\bm{\alpha}}$ is said to have degree $|\bm{\alpha}| = \alpha_1 + \alpha_2 + \cdots + \alpha_m$ and $\clh_k$ denotes the subspace in the Hilbert module $\clh$ spanned by the monomials having degree $k$ for $k$ in $\mathbb{N}$. 

If $(\bm{A})$ holds for $\clh$, then $\clh = \clh_0 \oplus \clh_1 \oplus \cdots $. Polynomials in $\clh_k$ are said to be {\it homogeneous} of degree $k$. An ideal $I$ in $\mathbb{C}[\mathbf{z}]$ is said to be a {\it homogeneous ideal} if it is generated by a set of homogeneous polynomials. 

For $I$ a homogeneous ideal in $\mathbb{C}[\mathbf{z}]$, let $[I]$ denote the closure of $I$ in the Hilbert module $\clh$. Then 
\begin{equation}\label{I}
 [I] = I_0 \oplus I_1 \oplus \cdots, 
\end{equation}
where $I_k = [I]\cap\clh_k$, and 
\begin{equation}\label{Iperp}
 [I]^{\perp} = I^{\perp}_0 \oplus I^{\perp}_1 \oplus \cdots,
\end{equation}
where $I^{\perp}_k = [I]^{\perp}\cap\clh_k$. 

Sometimes the following assumption, weaker than ($\bm{A}$), is sufficient to prove results:

\vspace{0.1in}

{\bf(A*)} the $\{\clh_k\}$ are orthogonal and $\clh = \clh_0 \oplus \clh_1 \oplus \cdots$. 

\vspace{0.1in}

The closure of a principal homogeneous ideal yields an example of a Hilbert modules satisfying ($\bm{A^*}$) in view of (\ref{I}). In fact, the same is true if $\clh$ only satisfies ($\bm{A^*}$). 

\vspace{0.2in}

\begin{Lemma}\label{lemmaA*}
 If $\clh$ is a Hilbert module satisfying ($\bm{A^*}$) and $I$ is a principal homogeneous ideal, then $[I]$ is a Hilbert module also  satisfying ($\bm{A^*}$).  
\end{Lemma}

\vspace{0.2in}

\begin{Remark} If $I$ is not principal, it is still finitely generates and $[I]$ satisfies a variant of ($\bm{A^*}$), where one now considers a Hilbert space completion of $\mathbb{C}[\bm{z}] \otimes \mathbb{C}^r$ for $r$ equal to the number of generators. We will offer only limited development of this idea in this paper.
\end{Remark}
\vspace{0.2in}

Another way to characterize homogeneous polynomials in a weighted shift Hilbert module $\clh$ is in terms of a natural unitary action of $\mathbb{T} = \mathbb{R}/2 \pi \mathbb{Z}$ on $\clh$. For $\lambda$ in $\mathbb{R}$ define $\gamma_{\lambda} (\mathbf{z}^{\bm{\alpha}}) = e^{i |{\bm{\alpha}}| \lambda} \mathbf{z}^{\bm{\alpha}}$ for $\bm{\alpha}$ in $\mathbb{N}^m$ or, equivalently, define $\gamma_{\lambda}(z_1, z_2, \ldots, z_m) = (e^{i \lambda} z_1, e^{i \lambda} z_2, \ldots, e^{i \lambda} z_m)$. Then $\clh_k$ is the eigenspace for this action for the eigen-character corresponding to $k$ in $\mathbb{Z} = \hat{\mathbb{T}}$. 

A generalization of the notion of homogeneous polynomial can be defined for an $m$-tuple of positive integers $(n_1, n_2, \ldots, n_m)$ or a weight $\bm{n}$. These polynomials can also be defined using the representation of $\mathbb{T} = \mathbb{R}/ {2 \pi \mathbb{Z}}$ on $\clh$ so that $$\gamma_{\lambda}^{\bm{n}} (z_1, z_2, \ldots, z_m) = (e^{i n_1 \lambda} z_1, e^{i n_2 \lambda} z_2, \ldots, e^{i n_m  \lambda} z_m).$$
The polynomials in the eigenspace $\clh_l^{\bm{n}}$ for this actions are the {\it quasi-homogeneous polynomials of degree $l$} for the weight $\bm{n}$.  Again, if $\clh$ satisfies ($\bm{A}$), we have that 
\begin{equation}
 \clh = \clh_0^{\bm{n}} \oplus \clh_1^{\bm{n}} \oplus \cdots
\end{equation}

\vspace{0.2in}

\begin{Remark}There is an obvious analogue of ($\bm{A^*}$) relevant for the consideration of quasi-homogeneous ideals. Given a weight $\bm{n}$, the Hilbert module $\clh$ over $\mathbb{C}[\bm{z}]$ satisfies
{\bf ($\bf{A_n^*}$)} if the subspaces $\{\clh^{\bm{n}}_l\}$ are orthogonal and $\clh = \clh^{\bm{n}}_0 \oplus \clh^{\bm{n}}_1 \oplus \cdots$ (Here we are defining the $\clh^{\bm{n}}_l$ directly using the weight $\bm{n}$ and not in terms of the circle group action.)
\end{Remark}

\vspace{0.2in}

Assume $\clh$ satisfies ($\bm{A}$). For a monomial $\mathbf{z}^{\bm{\alpha}}$, let $\clh^{\bm{n}}(\bm{\alpha})$ be the smallest subspace of $\clh$ containing $\mathbf{z}^{\bm{\alpha}}$ which is invariant under $M_{z_i^{n_i}}$ for $i = 1, 2, \ldots, m$. We want to make several observations about this family of subspaces of $\clh$. First, the subspaces $\clh^{\bm{n}}(\bm{\alpha})$ for $0 \leq \bm{\alpha} < \bm{n}$ are pairwise orthogonal and their direct sum is $\clh$, where $0 \leq \bm{\alpha} < \bm{n}$ means that $0 \leq \alpha_i < n_i$ for $i = 1, 2, \ldots, m$. Second, these subspaces are the minimal common reducing subspace of the $m$-tuple $\{M_{z_1}^{n_1}, M_{z_2}^{n_2}, \ldots, M_{z_m}^{n_m}\}$. Finally, we can make each $\clh^{\bm{n}}(\bm{\alpha})$ into the Hilbert module $\hat{\clh}^{\bm{n}}(\bm{\alpha})$ over $\mathbb{C}[\mathbf{z}]$ by defining$$ (p \cdot f) (\mathbf{z}) = p(z_1^{n_1}, z_2^{n_2}, \ldots, z_m^{n_m}) f(z_1, z_2, \ldots, z_m),$$ for $f$ in $\clh^{\bm{n}}(\bm{\alpha})$. We call this module action the {\it weighted module action determined by} $\bm{n}$. Moreover, $\hat{H}^{\bm{n}}(\bm{\alpha})$ is a Hilbert module over $\mathbb{C}[\bm{z}]$ which satisfies ($\bm{A}$) since $\clh$ does. We summerize this discussion in 

\vspace{0.2in}

\begin{Lemma}\label{gradedA}
 If $\clh$ satisfies ($\bm{A}$), then 

\begin{equation}\label{Hnalpha}
\clh = \oplus_{0 \leq \bm{\alpha} < \bm{n}}  \clh^{\bm{n}}(\bm{\alpha})
\end{equation}
 and each $\clh^{\bm{n}}(\bm{\alpha})$ is a Hilbert module $\hat{\clh}^{\bm{n}}(\bm{\alpha})$ over $\mathbb{C}[\mathbf{z}]$ for the weighted module action and $\hat{\clh}^{\bm{n}}(\bm{\alpha})$ satisfies ($\bm{A}$).
\end{Lemma}

\vspace{0.2in}

An ideal $J$ in $\mathbb{C}[\mathbf{z}]$ is said to be {\it quasi-homogeneous for the weight} $\bm{n}$ if it is generated by a set of quasi-homogeneous polynomials all for the weight $\bm{n}$; that is, the generators are in $\cup_{l =0}^{\infty} \clh_l^{\bm{n}}$ for some fixed weight $\bm{n}$. We call $\bm{n}$ the weight of the ideal. 

The analogue of Lemma \ref{lemmaA*} holds in that the closure of a principal quasi-homogeneous ideal $J$ of weight $\bm{n}$ will satisfy ($\bm{A_n^*}$).

Let $J$ be a quasi-homogeneous ideal of weight $\bm{n}$ with generators $q_1, q_2, \ldots, q_r$ in $\mathbb{C}[\mathbf{z}]$. For each $j = 1, 2, \ldots, r$, we have $<q_j> = \cup_{0 \leq \bm{\beta} \leq \bm{n} - \bm{\alpha}}  <{\mathbf{z}}^{\bm{\beta}} q_j>_{\bm{n}}$, where $<\; >$ denotes the ideal in $\mathbb{C}[\mathbf{z}]$ generated by the set in brackets and $< \;>_{\bm{n}}$ denotes the ideal in $\mathbb{C}[\mathbf{z}]$ with the weighted module action for weight $\bm{n}$ generated by the set in brackets. Moreover, we have $$J = <q_1, q_2, \ldots, q_r> = \vee_{j = 1}^{r} <q_{j}> = \vee_{j=1}^{r} \cup_{0 \leq \bm{\beta} \leq \bm{n} - \bm{\alpha}} < {\mathbf{z}}^{\bm{\beta}} q_j>_{\bm{n}} .$$

This identity represents $J$ as the finite direct sum of ideals $J_{\bm{\alpha}} = J \cap \clh^{\bm{n}}(\bm{\alpha})$ in $\clh^{\bm{n}}(\bm{\alpha})$, where the ideal structure is relative to the weighted module action for weight $\bm{n}$. 

A consequence of this decomposition is that some problems concerning quasi-homogeneous ideals can be reduced to the corresponding problems about homogeneous ideals. However, the results about the homogeneous case must be robust enough to cover the Hilbert modules obtained in the decomposition. 

\vspace{0.6in}

\newsection{Essentially reductive Hilbert modules}

\vspace{0.4in}

We are mainly interested in essentially reductive (or essentially normal) Hilbert modules; that is, 

\vspace{0.1in}

{\bf (B)} a Hilbert module $\clh$ over $\mathbb{C}[\mathbf{z}]$ such that the commutator $[M_f, M_g^*] = M_f M_g^* - M_g^* M_f$ is compact for $f, g$ in $\mathbb{C}[\mathbf{z}]$.

\vspace{0.1in}

A strengthening of this condition is 

\vspace{0.1in}

{\bf ($\bf{B_p}$)} a Hilbert module $\clh$ over $\mathbb{C}[\bm{z}]$ such that the commutator $[M_f, M_g^*]$ is in $\cll^p$ for some fixed $p$, $1\leq p \leq \infty$, and $\cll^p$ denotes the Schatten-von Neumann class. 

\vspace{0.1in}

Note that $(\bm{B_{\infty}}) = (\bm{B})$ and $(\bm{B_p})$ implies $(\bm{B})$. One knows that the Hardy, Bergman and $m$-shift Hilbert modules satisfy ($\bm{B_p}$) for $p >m$. Arveson raised the question of whether the closure $[I]$ of a homogeneous ideal $I$ in $H_m^2$, satisfies ($\bm{B_p}$), for $p>m$. Guo and Wang extended the question to include the closure of quasi-homogeneous ideals. We show that the second question can be reduced to the first one if  the affirmative solution for the homogeneous case is robust enough.

We next want to place more restrictions on $\clh$ so that it resembles more closely the Hardy, Bergman and $m$-shift Hilbert  modules and relates directly to the questions raised in the introduction.

Assume that ($\bm{B}$) holds so that $\clh$ is essentially reductive. Then the $C^*$-algebra $\clt(\clh)$ generated by the operators $\{M_p : p(\mathbf{z}) \in \mathbb{C}[\mathbf{z}]\}$  is a $C^*$-extension of the subalgebra $\clk(\clh)$ of compact operators on $\clh$ by $C(X)$. Here $X$ is some compact metrizable space which can be identified as a subset of $\mathbb{C}^m$; that is, $\clt(\clh)/\clk(\clh) \cong C(X)$ (see \cite{BDF}). Note that $\clt(\clh)$ contains a non-zero compact operator since $[M_{z_i}, M_{z_i}^*]$ is compact and non-zero. Hence, $\clk(\clh)$ is contained in $\clt(\clh)$ since $\clh$ is irreducible (cf. \cite{DB}). 

We make the further assumption that 

\vspace{0.1in}

{\bf (C*)} $\clh$ satisfies ($\bm{B}$) and $(\overline{M}_{z_1}, \overline{M}_{z_2}, \ldots, \overline{M}_{z_m})$ identifies $X$ as a subset of the unit sphere $\partial \mathbb{B}^m$, where $\overline{M}_{z_i}$ denotes the image of ${M}_{z_1}$ in the quotient algebra $\clt(\clh)/\clk(\clh)$.

\vspace{0.1in}

One can also assume 

\vspace{0.1in}

{\bf (C**)} $I_{\clh} - \sum_{i=1}^{m} M_{z_i}^* M_{z_i}$ is compact or $(M_{z_1}, M_{z_2}, \ldots, M_{z_m})$ is an essential spherical isometry. 

\vspace{0.1in}

In the presence of ($\bm{B}$), ($\bm{C^*}$) is equivalent to ($\bm{C^{**}}$). We also consider the following apparent strengthening of ($\bm{C^*}$).

\vspace{0.1in}

{\bf (C)} $\clh$ is a Hilbert module satisfying ($\bm{B}$) such that  $\clt(\clh)/\clk(\clh) \cong C(\partial \mathbb{B}^m)$.

\vspace{0.1in}

Although we will make little use of the following  notion in this paper, we include it for completeness:

\vspace{0.1in}
{\bf ($\bf{C_p}$)}  $\clh$ satisfies ($\bm{B_p}$) and $I_{\clh} - \sum_{i=1}^{m} M_{z_i}^* M_{z_i}$ is in $\cll^p$ for $1 \leq p \leq \infty$. 

\vspace{0.1in}

We show that for a weighted shift Hilbert module ($\bm{C^*}$) implies ($\bm{C}$) after establishing the following lemma.

\vspace{0.2in}

\begin{Lemma}\label{indexAB}
 If $\clh$ is a Hilbert module satisfying ($\bm{A^*}$) and ($\bm{B}$), so that the $m$-tuple $(M_{z_1}, M_{z_2}, \ldots\\, M_{z_m})$ is Fredholm, then the index of the corresponding Koszul complex is $-1$. Moreover, if $\clh$ also satisfies ($\bm{C^*}$), then it satisfies ($\bm{C}$).
\end{Lemma}

\vspace{0.2in}

\NI {\sf Proof.} Using the orthogonal direct sum decomposition $\clh = \clh_0 \oplus \clh_1 \oplus \cdots$, the Koszul complex for $(M_{z_1}, M_{z_2}, \ldots, M_{z_m})$  can be reduced to the direct sum of the corresponding Koszul complexes for the action of $(z_1, z_2, \ldots, z_m)$ on $\mathbb{C}[\mathbf{z}]$ which has index $-1$. The key to this reduction depends on the Fredholmness assumption of the complex which implies that all the maps in the Koszul complex for $\clh$ have closed range. Hence, the existence of approximate solutions implies that there exists a solution. Since the Fredholm index of the $m$-tuple $(M_{z_1}, M_{z_2}, \ldots, M_{z_m})$ is not zero, the $C^*$-extension of $\clk(\clh)$ defined by $\clt(\clh)$ is nontrivial. Hence $X$ can not be a proper subset of $\partial \mathbb{B}^m$, since the map from $K_1(X)$ to $K_1(\partial \mathbb{B}^m)$ is the zero map for such $X$. \qed

\vspace{0.3in}

Note that a related result appears in \cite{GRS} but with a different argument. 

Suppose $\clh$ is a Hilbert module over $\mathbb{C}[\bm{z}]$ such that the coordinate multiplier $m$-tuple $(M_{z_1}, M_{z_2}, \ldots, M_{z_m})$ yields an essentially spherical isometry. The question of whether or not this $m$-tuple is an essentially spherical unitary is equivalent to the question of whether $\clh$ satisfies ($\bm{B}$) in view of the fact that 
\begin{equation}\label{M_i}
I - \sum_{i=1}^{m} [\bar{M_{z_i}}, \bar{M_{z_i}}^*] = \sum_{i=1}^{m} \bar{M_{z_i}}^* \bar{M_{z_i}} - \sum_{i=1}^{m} \bar{M_{z_i}} \bar{M_{z_i}}^* = \sum_{i=1}^{m} [\bar{M_{z_i}}^*, \bar{M_{z_i}}],
\end{equation}

\NI where $\bar{M_{z_i}}$ denotes the operator $M_{z_i}$ modulo the compacts. Here we are using the fact that the operators $\{M_{z_1}, M_{z_2}, \ldots, M_{z_m} \}$ are essentially hyponormal. 

To relate the results for quasi-homogeneous ideals to those for homogeneous ideals, we need to relate condition ($\bm{C}$) for the action of $\mathbb{C}[\mathbf{z}]$ on $\clh$ to the corresponding condition for the weighted module action of $\mathbb{C}[\bm{z}]$ on the $\clh^{\bm{n}} (\bm{\alpha})$ for weight $\bm{n}$ and $\bm{\alpha}$ in $\mathbb{N}^m$.

\vspace{0.2in}

\begin{Lemma}
 A Hilbert module satisfying ($\bm{B_p}$) satisfies ($\bm{C_p}$) if and only if 

\vspace{0.1in}

{$\bm{(D_p)}$} $(M_{z_1}, M_{z_2}, \ldots, M_{z_m})$ is an essentially spherical isometry and $M_{z_1}^* M_{z_1} + M_{z_2}^* M_{z_2}+ \cdots + M_{z_m}^* M_{z_m} - I_{\clh}$ is in $\cll^p$ where $1\leq p \leq \infty$.
\end{Lemma}

\vspace{0.2in}

\NI {\sf Proof.} In $\clt(\clh)/\clk(\clh)$ the image of $M_{z_1}^* M_{z_1} + M_{z_2}^* M_{z_2}+ \cdots + M_{z_m}^* M_{z_m} - I_{\clh}$ is $\sum_{i=1}^{m} |z_i|^2 - 1$, which vanishes on $\partial \mathbb{B}^m$. Identity (\ref{M_i}) completes the proof. \qed

\vspace{0.3in}

Not all Hilbert modules satisfying ($\bm{A}$) and ($\bm{B}$) also satisfy ($\bm{C}$). Consider the Hilbert modules based on Reinhardt domains which, in general, do not satisfy ($\bm{C}$). In particular, the maximal ideal spaces in these cases are not always $\partial \mathbb{B}^m$. We will discuss later how the relationship between ($\bm{B_p}$) and ($\bm{D_p}$) relates to the conjecture of Arveson and the refinement of it by the first author. Now we want to continue developing the relation between quasi-homogeneous ideals and related homogeneous ones which we began in Lemma \ref{gradedA}.

\vspace{0.2in}

\begin{Lemma}\label{gradedBC}
 Suppose $\clh$ satisfies ($\bm{A}$), where $\bm{n}$ is a weight, and $$\clh = \oplus_{0 \leq \bm{\alpha} < \bm{n}}  \clh^{\bm{n}}(\bm{\alpha})$$ is the decomposition in (\ref{Hnalpha}). Then $\clh$ satisfies ($\bm{B_p}$) or ($\bm{C_p}$), if and only if all the $\hat{\clh}^{\bm{n}}(\bm{\alpha})$, for $0 \leq \bm{\alpha} <\bm{n}$, satisfy ($\bm{B_p}$) and ($\bm{C_p}$), where $1\leq p \leq \infty$. 
\end{Lemma}

\vspace{0.2in}

\NI {\textsf Proof.} Since each of the $\clh^{\bm{n}}(\bm{\alpha})$ reduces the $m$-tuple $(M_{z_1}^{n_1}, M_{z_2}^{n_2}, \ldots, M_{z_m}^{n_m})$, we can express the commutators as an orthogonal direct sum, 
\begin{equation}\label{Mcap1}
 [M_{z_i^{n_i}}, M_{z_j^{n_j}}^*] = \oplus_{0 \leq \bm{\alpha} < \bm{n}} [\hat{M}_{z_i}, \hat{M}_{z_j}^*],
\end{equation}
where $\hat{M}_{z_i}$ is defined by the weighted module action of $z_i$ on $\hat{\clh}^{\bm{n}}(\bm{\alpha})$. Similarly, we have
\begin{equation}\label{Mcap2}
 \sum_{i=1}^{m} M_{z_i^{n_i}}^* M_{z_i^{n_i}} =\oplus_{0 \leq \bm{\alpha} < \bm{n}} \{\sum_{i=1}^{m} \hat{M_{z_i}}^* \hat{M_{z_i}}\}.
\end{equation}
Thus, from (\ref{Mcap1}) and (\ref{Mcap2}), we see that $\clh$ satisfies ($\bm{B_p}$) or ($\bm{C^*_p}$) if and only if all the $\hat{\clh}^{\bm{n}}(\bm{\alpha})$ do for $0 \leq \bm{\alpha} < \bm{n}$. Thus the Koszul complex for the $m$-tuple is exact, or Fredholm, if and only if the same is true for each of the Koszul complexes for the direct summands. Thus $\clh$ satisfying ($\bm{C_p}$) implies that each $\hat{\clh}^{\bm{n}}(\bm{\alpha})$ satisfies ($\bm{C^*_p}$). But Lemma \ref{indexAB} implies that each $\hat{\clh}^{\bm{n}}(\bm{\alpha})$ satisfies ($\bm{C}$), which concludes the proof. \qed

\vspace{0.3in}

\begin{Remark}If $\clh$ satisfies ($\bm{A_n^*}$), then the $\hat{\clh}^{\bm{n}}(\bm{\alpha})$ will satisfy ($\bm{A^*}$) and the ideals $J^{\bm{n}}(\bm{\alpha}) = J \cap \hat{\clh}^{\bm{n}}(\bm{\alpha})$ are homogeneous in $\hat{\clh}^{\bm{n}}(\bm{\alpha})$ relative to the weighted module action.
\end{Remark}

\vspace{0.2in}

We need a lemma concerning hyponormal operators to complete the reduction of questions concerning quasi-homogeneous ideals to the  analogous questions about homogeneous ones.

\vspace{0.2in}

\begin{Lemma}
 If $T$ is an essentially hyponormal operator for which $T^k$ is essentially normal for some $k \geq 1$, then $T$ is essentially normal.
\end{Lemma}
\vspace{0.1in}

\NI {\textsf Proof.} Working modulo the compacts in the Calkin algebra, the question reduces to showing that a hyponormal operator $\bar{T}$ for which $\bar{T}^k$ is normal must itself be normal. To that end, consider the spectral representation for $\bar{T}^k$ so that $$\bar{T}^k = \int_{\sigma(\bar{T}^k)} z \; \mbox{d}E_z,$$

\NI for a spectral measure $\{E_z\}$ on $\sigma(\bar{T}^k)$. Since $\bar{T}$ commutes with $\bar{T}^k$, there exists a measurable operator-valued function $X(z)$ so that $\bar{T} = \int_{\sigma(\bar{T}^k)} X(z) \; \mbox{d}E_z$. But $\bar{T}$ is hyponormal if and only if $X(z)$ is hyponormal a.e. and $X(z)^k = z I$ a.e. This implies that $X(z)$ is a normal operator a.e. with spectrum contained in the set of $k$-th roots of $z$, which implies that $\bar{T} =  \int_{\sigma(\bar{T}^k)} X(z) \; \mbox{d}E_z$ is normal. \qed

\vspace{0.2in}

Collecting the lemmas, we have the following reduction.

\vspace{0.2in}

\begin{Theorem}\label{BC}
 Assume that $\clh$ is a Hilbert module satisfying ($\bm{A}$), that each $M_{z_i}, \; i = 1, 2, \ldots, m$ is essentially hyponormal and that $J$ is a quasi-homogeneous ideal for weight $\bm{n}$. Consider the decomposition  

$$\clh = \oplus_{0 \leq \bm{\alpha} < \bm{n}}  \clh^{\bm{n}}(\bm{\alpha})$$ 

\NI and set $ J^{\bm{n}}(\bm{\alpha}) = J \cap  \clh^{\bm{n}}(\bm{\alpha})$ for $0 \leq \bm{\alpha} < \bm{n}$. Then $[J] = \oplus_{0 \leq \bm{\alpha} < \bm{n}} [J^{\bm{n}} (\bm{\alpha}) ]$ with each $[J^{\bm{n}} (\bm{\alpha})]$ being a homogeneous ideal for $\mathbb{C}[\mathbf{z}]$ with the weighted module action. Then $[J]$ satisfies ($\bm{B_p}$) or ($\bm{B_p}$) and ($\bm{C_p}$) if and only if each $[J^{\bm{n}} (\bm{\alpha})]$ does for $0 \leq \bm{\alpha} < \bm{n}$, with $1 \leq p \leq \infty$.
\end{Theorem}

\vspace{0.1in}
\NI {\textsf Proof.} The earlier lemmas yield the first decomposition and now for $J$ as well because of the relation between cross-commutators $[M_{z_i^{n_i}}, M_{z_j^{n_j}}^*]$ on each $[J^{\bm{n}} (\bm{\alpha})]$ and $[M_{z_i^{n_i}}, M_{z_j^{n_j}} ^*]$ on $J$. Therefore, one sees that the restriction of the operators $\{M_{z_i^{n_i}}\}$ to $[J]$ have $\cll^p$ commutators. Now the fact that each $M_{z_i}$ is essentially hyponormal implies that $M_{z_i}|_{[J]}$ is essentially hyponormal and the previous lemma completes the proof. \qed

\vspace{0.3in}

\begin{Corollary}\label{BCcor}
 If $\clh$ is an isometric Hilbert module over $\mathbb{C}[\bm{z}]$ satisfying ($\bm{A_n^*}$) and $J$ is a quasi-homogeneous ideal for weight $\bm{n}$, then $[J]$ satisfies ($\bm{B_p}$) or ($\bm{B_p}$) and ($\bm{C_p}$) if and only if each $[J^{\bm{n}}(\bm{\alpha})]$ does for $0 \leq \bm{\alpha} < \bm{n}$ and $1 \leq p \leq \infty$.
\end{Corollary}
\vspace{0.1in}

\NI \textsf{Proof.} In case $\clh$ is actually isometric, then it is subnormal by \cite{A} and hence the $M_{z_i}$ are hyponormal for $i=1, 2, \ldots, m$. Thus the Theorem applies. \qed

\vspace{0.3in}

One is tempted to conclude that this result allows one to apply the recent results of Guo and Wang \cite{GW} on homogeneous ideals to quasi-homogeneous ones. However, those results are particular to the $m$-shift space. Unfortunately, the Hilbert modules in the decomposition given in Lemma \ref{gradedA} do not have the unitary symmetry necessary to allow one to apply these technique.

In this section we have considered the case of multiplicity one; that is, Hilbert module, obtained as the completion of $\mathbb{C}[\mathbf{z}]$. It would be natural to consider completions of $\mathbb{C}[\mathbf{z}] \otimes \mathbb{C}^k$ for $k$ in $\mathbb{N}$ as mentioned in Section 1. We will return to this issue in Section 5 in order to state stronger theorems.

\vspace{0.6in}

\newsection{The case $m=2$}

\vspace{0.4in}

Our goal in this section is to prove the following theorem for the case of $m=2$ and multiplicity one.

\vspace{0.2in}

\begin{Theorem}\label{BCD}
Let $\clh$ be a Hilbert module completion of $\mathbb{C}[z_1, z_2]$ which satisfies ($\bm{A^*}$), ($\bm{B}$) and ($\bm{C}$) and $I$ be a homogeneous ideal in $\mathbb{C}[z_1, z_2]$. Then $[I]$ and $\clh/[I]$ are essentially reductive or, equivalently, $[I]$ satisfies ($\bm{D}$).
\end{Theorem}

\vspace{0.3in}

In case $\clh$ is $H^2_2$, this result was proved by Guo in \cite{G}. Moreover, the result is known to hold for $\clh = H^2_2 \otimes \mathbb{C}^r$ and closely related Hilbert modules. This was established by Guo and Wang (see \cite{GW}). However, the techniques in those papers do not seem to extend to yield the result for Hilbert modules as genearl as these considered here.. 

\vspace{0.2in}

\begin{Theorem}\label{ABCDtheorem}
Let $\clh$ be a Hilbert module completion of $\mathbb{C}[z_1, z_2]$ which satisfies ($\bm{A}$), ($\bm{B}$) and ($\bm{C}$) and $J$ be a quasi-homogeneous ideal in $\mathbb{C}[z_1, z_2]$ having weight $\bm{n} = (n_1, n_2)$. Then $[J]$ and $\clh/[J]$ are essentially reductive and satisfy ($\bm{D}$).
\end{Theorem}
\vspace{0.1in}

\NI \textsf{Proof.} The result follows from Theorem \ref{BCD} in view of Theorem \ref{BC} and Lemma \ref{gradedBC}. \qed

\vspace{0.3in}

This result was obtained by Guo and Wang in \cite{GW} in the case $\clh = H^2_2$. To prove Theorem \ref{BCD}, we adopt the outline of their proof but the necessary lemmas are established by different means.

\vspace{0.3in}
\begin{Lemma}
If $\clh$ is a Hilbert module completion of $\mathbb{C}[z_1, z_2]$ satisfying ($\bm{B}$) and ($\bm{C}$) then $[z_1 \clh]$ satisfies ($\bm{B}$) and ($\bm{C}$).
\end{Lemma}

\NI \textsf{Proof.} Note that $[z_1\clh]^{\perp} = \mbox{ker} M_{z_1}^*$. If we decompose $M^*_{z_1}$ relative to $\clh = \mbox{ker} M^*_{z_1} \oplus [z_1 \clh]$, then we obtain the matrix 
\[ \left( \begin{array}{cc}
0 & A \\
0 & B \end{array} \right).\] 

Since $M^*_{z_1}$ is essentially normal, it follows that $A$ is compact and $B$ is essentially normal. If we write $M_{z_2}^*$ as 

\[ \left( \begin{array}{cc}
C & D \\
0 & E \end{array} \right),\] 

\NI then the fact that $I_{\clh} - M^*_{z_1} M_{z_1} -  M_{z_2}^* M_{z_2}$ is compact implies that $I - B B^* - E E^*$ is compact. Moreover, the fact that $[M_{z_1}^*, M_{z_2}^*] = 0$ implies that $B E = E B$ and the fact that $B$ is essentially normal implies that $B^* E = E B^*$. In conclusion we have that $E$ and $E^*$ commutes modulo the compacts with $I - B B^* = E E^*$ which implies via a polar form argument that $E$ is essentially normal. This completes the proof. \qed

\vspace{0.2in}

\begin{Remark}
One can generalize the preceding proof to show that $[<z_1, z_2, \ldots, z_{m-1}>]$ is essentially reductive where the closure is taken in a Hilbert module $\clh$ over $\mathbb{C}[z_1, z_2, \ldots, z_m]$ satisfying ($\bm{B}$) and ($\bm{C}$). 
\end{Remark}

\vspace{0.3in}

\begin{Lemma}\label{zalpha}
Let $\clh$ be a Hilbert module completion of $\mathbb{C}[z_1, z_2]$ satisfying ($\bm{A^*}$), ($\bm{B}$) and ($\bm{C}$). The ideal $I = <z_1 + \alpha z_2>$  generates the submodule $[I]$ of $\clh$ and quotient module $\clh/[I]$, both of which satisfy ($\bm{A^*}$), ($\bm{B}$) and ($\bm{C}$).
\end{Lemma}
\vspace{0.1in}

\NI \textsf{Proof.} We observe first that the preceding lemma handle the case $\alpha = 0$. Next, we observe that ($\bm{A^*}$) holds for $[I]$ using Lemma \ref{lemmaA*} and the fact that $I$ is principal. Once we establish ($\bm{B}$) for $[I]$, we will see that ($\bm{C}$) follows. Let $(\overline{M}_{z_1}, \overline{M}_{z_2})$ denote the coordinate multiplication operators defined on the quotient module $\clh/[I]$. Suppose $\alpha \neq 0$. Since $\clh$ satisfies ($\bm{B}$) and ($\bm{C}$), we have that $M_{z_1} M_{z_1}^* + M_{z_2} M_{z_2}^* -I_{\clh}$ is compact or that $(M_{z_1}^*, M_{z_2}^*)$ is an essentially spherical isometry. A simple matrix computation shows that the restriction of $(M_{z_1}^*, M_{z_2}^*)$ to $\mathcal{Q} = \clh/[I] = \clh \ominus [I]$ is also an essentially spherical isometry or $\overline{M}_{z_1} \overline{M}_{z_1}^* + \overline{M}_{z_2} \overline{M}_{z_2}^* -I_{\mathcal{Q}}$ is compact. But since $z_1 + \alpha z_2$ is in $I$, we see that $\overline{M}_{z_1} + \alpha \overline{M}_{z_2} = 0$ or $\overline{M}_{z_2} = - \frac{1}{\alpha} \overline{M}_{z_1}$ and therefore, $\overline{M}_{z_1} \overline{M}_{z_1}^* + \frac{1}{| \alpha |^2} \overline{M}_{z_2} \overline{M}_{z_2}^* -I_{\mathcal{Q}}$ is compact. Hence, $(1 + |\alpha|^2)^{- \frac{1}{2}} \bar{M}_{z_1}^*$ is an essential isometry. Similarly, we have that $(1 + |\alpha|^2)^{- \frac{1}{2}} \overline{M}_{z_2}$ is an essential isometry. Easy calculation shows that $[I]_k^{\perp}$ is one-dimensional and that $\bar{M}_{z_1}$ on $\clh/[I]$ is a weighted shift (cf. \cite{DMcal}). The fact that $(1+ |\alpha|^2)^{-\frac{1}{2}} \bar{M}^*_{z_1}$ is an essential isometry shows that the absolute value of the weights converge to 1. Hence $\bar{M}^*_{z_1}$ is essentially normal which implies that $[I]$ and $\clh/[I]$ satisfy ($\bm{B}$). This completes the proof of the lemma. \qed

\vspace{0.3in}

For the proof of Theorem \ref{BCD}, we need the following lemma.

\vspace{0.2in}

\begin{Lemma}
Suppose $\clh$ is a Hilbert module completion of $\mathbb{C}[\bm{z}]$ satisfying ($\bm{B}$) and ($\bm{C}$). If $I = <p(\bm{z})>$ is a principal ideal satisfying ($\bm{B}$), then $[I]$ also satisfies ($\bm{C}$).
\end{Lemma}
\vspace{0.1in}
\NI \textsf{Proof.} Observe first that one can identify $$\clt(\clh)/\clk(\clh) \cong \clt([I])/\clk([I]) \oplus \clt(\clh/[I])/\clk(\clh/[I])$$ since the off-diagonal entries for the matrix representation for the operators in $\clt(\clh)$ are compact. The union of the two maximal ideal spaces equals $\partial \mathbb{B}^m$. Further, the maximal ideal space for the quotient algebra is the intersection of $\partial \mathbb{B}^m$ with the zero variety $\clz$ of $p(\bm{z})$ (see \cite{GRS}). Therefore, we see that $\clt([I])/\clk([I])$ satisfies ($\bm{C^*}$). Finally, the image of the odd $K$-homology element defined by the quotient module is 0 because the maximal ideal space is a proper subset of $\partial \mathbb{B}^m$. Thus we see that $[I]$ satisfies ($\bm{C}$). \qed

\vspace{0.3in}
 One advantage in working with general Hilbert modules, rather than a specific one such as the $m$-shift space, is that induction can be used as follows:

\vspace{0.2in}

\begin{Lemma}
Assume for all Hilbert module completions of $\mathbb{C}[\bm{z}]$ satisfying ($\bm{A^*}$), ($\bm{B}$) and ($\bm{C}$), that the principal ideals $I_1 = <p_1(\bm{z})>$ and $I_2 = <p_2(\bm{z})>$ also satisfy ($\bm{A^*}$), ($\bm{B}$) and ($\bm{C}$). Then the same is true for the principal ideal $K = <p(\bm{z}) q(\bm{z})>$. 
\end{Lemma}
\vspace{0.1in}

\NI \textsf{Proof.} The proof follows once one observes that by assumption the closure $[I_1]$ of $I_1$ satisfies ($\bm{A^*}$), ($\bm{B}$) and ($\bm{C}$). Therefore, the closure of the ideal $I_2$ in $[I_1]$ satisfies ($\bm{A^*}$), ($\bm{B}$) and ($\bm{C}$), also by assumption (Here, we are using the weakly wandering vector $p(\bm{z})$ in $[I_1]$ to identify $[I_1]$ as a completion of $\mathbb{C}[\bm{z}]$.) But this closure equals $[<p(\bm(z) q(\bm{z})>]$ or the closure of $K$ in $\clh$, which completes the proof. \qed

\vspace{0.3in}

\NI \textbf{ Proof of Theorem 4.1:} Recall the fact (\cite{GW2}) that every homogeneous polynomial in $\mathbb{C}[z_1, z_2]$ is the product of a monomial $z_1^{k_1} z_2^{k_2}$ and factors of the form $(z_1 + \alpha z_2)$ with $\alpha \neq 0$. Since $\clh$ satisfies ($\bm{A^*}$), from Lemma \ref{zalpha}, we see that the closure of any homogeneous principal ideal in $\mathbb{C}[z_1, z_2]$ satisfies ($\bm{A^*}$), ($\bm{B}$) and ($\bm{C}$). The proof is completed by appealing to a result of Yang \cite{Y} showing that any proper homogeneous ideal in $\mathbb{C}[z_1, z_2]$ contains a homogeneous principal ideal of finite co-dimension which complete the proof. \qed

\vspace{0.6in}

\newsection{Quotient modules of dimension one}

\vspace{0.4in}

If $I$ is a homogeneous ideal in $\mathbb{C}[\bm{z}]$, then there is an intimate relation between the zero variety $\clz = \{ \bm{z} \in \mathbb{C}^m : p(\bm{z}) = 0 \;\; \mbox{for all} \;p \;\; \mbox{in} \; I\}$ of $I$ and the Hilbert-Samuel polynomial $p_{\mathbb{C}[\bm{z}]/I}$ for the quotient module $\mathbb{C}[\bm{z}]/I$. The same is true for Hilbert modules. In particular, $p_{\mathbb{C}[\bm{z}]/I}$ will be linear if and only if $\clz$ has complex dimension one. Rather than developing these ideas here, we will use a consequence of this fact as our basic assumption since we want to consider the case of higher multiplicity any way. Here, the notion of zero variety is more complex.

Let $\clh$ be a Hilbert module completion of $\mathbb{C}[\bm{z}] \otimes \mathbb{C}^r$, for $r \geq 1$. The degree of a monomial $\bm{z}^{\bm{\alpha}} \otimes u$ for $u$ in $\mathbb{C}^r$ is $|\bm{\alpha}|$. An element of $\mathbb{C}[\bm{z}] \otimes \mathbb{C}^r$ is said to be {\it homogeneous} if all monomials in it have the same degree. A submodule $\cls$ of $\clh$ is said to be {\it homogeneous} if it is generated by homogeneous elements of  $\mathbb{C}[\bm{z}] \otimes \mathbb{C}^r$.  One knows that $\cls \cap \mathbb{C}[\bm{z}] \otimes \mathbb{C}^r$ is finitely generated and since its closure is $\cls$, hence so is $\cls$. Here one is relying on homogeneity being characterized by the circle group action. One can extend ($\bm{A}$) and ($\bm{A^*}$) to such modules in an obvious fashion and again one has the orthogonal decompositions $\clh = \oplus \clh_k$, $\cls = \oplus \cls_k$ with $\cls_k = \cls \cap \clh_k$, and $\cls^{\perp} = \oplus \cls_k^{\perp}$ with $\cls_k^{\perp} = \cls^{\perp} \cap \clh_k$ as in the case of $r = 1$. 

For each $k$ in $\mathbb{N}$ and $i=1, 2, \ldots, m$, there exists an operator $A_{i,k} : \cls_k^{\perp} \raro \cls_{k+1}^{\perp}$, so that 

\begin{equation}\label{Ashift}
\overline{M}_{z_i} (p_0, p_1, p_2, \ldots) = (0, A_{i,0}p_0, A_{i,1}p_1, A_{i,2}p_2, \ldots),
\end{equation}
where $p_k$ is in $\cls^{\perp}_k$ and $\overline{M}_{z_i}$ is the compression of $M_{z_i}$ to the quotient module $\cls^{\perp}$. 

We say that $\cls^{\perp}$ has {\it bounded dimension} if $\mbox{dim} \cls^{\perp}_k \leq M < \infty$ for some $M$. Using the existence of the Hilbert-Samuel polynomial \cite{D-Y}, one observes that the assumption that dimension of $\cls^{\perp}$ is bounded implies the existence of natural numbers $M_0$ and $K$ such that dim$\cls_k^{\perp} = M_0$ for $k \geq K$. 

If $\cls^{\perp}$ has bounded dimension and $\clh$ satisfies ($\bm{A^*}$), then the operators $M_{z_i}$ are unilateral block weighted shifts for $i=1, 2, \ldots, m$. In particular, for all $k \geq K$, the operators $A_{i, k} : \cls^{\perp}_k \raro  \cls^{\perp}_{k+1}$ map between spaces of dimension $M_0$. Moreover, the adjoint $\overline{M}_{z_i}^*$ satisfies 

\begin{equation}\label{MandA}
\overline{M}_{z_i}^* (\oplus_{k=K}^{\infty} p_{k+1}) = \oplus_{k=K}^{\infty} A^*_{i, k+1} p_{k+1},
\end{equation}
\NI for $\oplus_{k=K}^{\infty} p_{k+1}$ in $\oplus_{k=K}^{\infty} \cls^{\perp}_{k+1}$, where $A^*_{i, k+1} : \cls^{\perp}_{k+1} \raro  \cls^{\perp}_{k}$.

For this calculation to be valid, it is essential that the $\cls_k^{\perp}$ be orthogonal. We begin with the following lemma.

\vspace{0.2in}

\begin{Lemma}
Let $\clh$ be a Hilbert module completion of $\mathbb{C}[\bm{z}] \otimes \mathbb{C}^r$ which satisfies ($\bm{B_p}$) and ($\bm{C_p^*}$) for some $1 \leq p \leq \infty$. Then the $m$-tuple $(M_{z_1}^*, M_{z_2}^*, \ldots, M_{z_m}^*)$ is an essentially spherical isometry and, moreover, $I - \sum_{i=1}^{m} M_{z_i} M_{z_i}^*$ is in $\cll^p$.
\end{Lemma}
\vspace{0.1in}
\NI \textsf{Proof.} The result follows from the identity 
\begin{equation}I - \sum_{i=1}^{m} M_{z_i} M_{z_i}^* = I - \sum_{i=1}^{m} M_{z_i}^* M_{z_i} + \sum_{i=1}^{m} [M_{z_i}^*, M_{z_i}],
\end{equation}
which we used in Section 3. \qed

\vspace{0.3in}

\begin{Lemma}\label{PN}
Let $\clh$ be a Hilbert module completion of $\mathbb{C}[\bm{z}] \otimes \mathbb{C}^r$ which satisfies ($\bm{B_p}$) and ($\bm{C_p^*}$) for some $1 \leq p \leq \infty$, $\cls$ be a homogeneous submodule of $\clh$ and $\overline{M}_{z_i}$ be the compression of $M_{z_i}$ to $\cls^{\perp}$ for $i=1, 2, \ldots, m$. Then $(\overline{M}_{z_1}, \overline{M}_{z_2}, \ldots, \overline{M}_{z_m})$ is an essentially spherical isometry and $I_{\cls^{\perp}} - \sum_{i=1}^{m} \overline{M}_{z_i} \overline{M}_{z_i}^*$ is in $\cll^p$. Moreover, for $i=1, 2, \ldots, m$, $[\overline{M}_{z_i}, \overline{M}_{z_i}^*] = P_i - N_i$, where $P_i$ and $N_i$ are positive and $N_i$ is $\cll^p$.
\end{Lemma}
\vspace{0.1in}

\NI \textsf{Proof.} The result follows by a matrix calculation using the fact that $I_{\clh} - \sum_{i=1}^{m} M_{z_i} M_{z_i}^*$ is in $\cll^p$ and $\cls^{\perp}$ is a joint invariant submodule for the $m$-tuple  $(M_{z_1}^*, M_{z_2}^*, \ldots, M_{z_m}^*)$. 
The last statement follows from the matrix calculation $[\bar{M}_{z_i}, \bar{M}_{z_i}^*] = P_{\cls^{\perp}} [{M}_{z_i}, {M}_{z_i}^*] P_{\cls^{\perp}} + (P_{\cls} M_{z_i} P_{\cls^{\perp}})^* (P_{\cls} M_{z_i} P_{\cls^{\perp}})$ since $[{M}_{z_i}, {M}_{z_i}^*]$ is in $\cll^p$. \qed

\vspace{0.3in}

\begin{Theorem}
Let $\clh$ be a Hilbert module completion of $\mathbb{C}[\bm{z}] \otimes \mathbb{C}^r$ satisfying ($\bm{A^*}$), ($\bm{B}$) and ($\bm{C}$) and $\cls$ be a homogeneous submodule of bounded dimension. Then $\cls$ and $\cls^{\perp}$ are essentially reductive.
\end{Theorem}
\vspace{0.1in}
\NI \textsf{Proof.} We begin by calculating the various operators involved using the representations in equations (\ref{Ashift}) and (\ref{MandA}). We consider only elements in $\cls_K^{\perp} \oplus \cls_{K+1}^{\perp} \oplus \cls_{K+2}^{\perp} \ldots$, which is sufficient since $\cls_0^{\perp} \oplus \cls_{1}^{\perp} \oplus \cdots \oplus \cls_{K-1}^{\perp}$ has finite dimension. In particular, we have $$\overline{M}_{z_i} (\oplus_{k = K}^{\infty}  p_k ) = \oplus_{k = K+1}^{\infty} A_{i, k} p_k,$$ and $$\overline{M}_{z_i}^* (\oplus_{k=K}^{\infty} p_{k+1}) = \oplus_{k=K}^{\infty} A^*_{i, k+1} p_{k+1},$$
where $p_k$ is in $\cls^{\perp}_k$. If we set $X = I - \sum_{i=1}^{m} \overline{M}_{z_i} \overline{M}_{z_i}^*$, then $X = \oplus X_k$, where $X_k$ is in $\cll(\cls_k^{\perp})$. Moreover, $\mbox{lim}_{k \raro  \infty} \|X_k\| = 0$, since $X$ is compact. 
Further, if we write $[\overline{M}_{z_i}, \overline{M}_{z_i}^*] = P_i - N_i$ as in Lemma \ref{PN}, we have $P_i = \oplus P_{i,k}$, and $N_i = \oplus N_{i, k}$ for $P_{i,k}, N_{i,k}$ in $\cll(\cls_k^{\perp})$ for $k$ in $\mathbb{N}$. Therefore, we have 

\begin{align*}
I_{\cls^{\perp}} - \sum_{i=1}^{m} \overline{M}_{z_i}^* \overline{M}_{z_i} & = \oplus (I_{\cls^{\perp}_k} - \sum_{i=1}^{m} A_{i, k}^* A_{i, k}) \\ & = \oplus (I_{\cls^{\perp}_k} - \sum_{i=1}^{m} A_{i, k-1} A_{i, k-1}^*) + \oplus (A_{i,k}^* A_{i,k} - A_{i,k-1} A_{i,k-1}^*)  \\ & = \oplus X_{k} + \oplus \sum_{i=1}^{m}(P_{i, k} - N_{i, k}) = \oplus ( X_{k} + \sum_{i=1}^{m} P_{i, k} - N_{i, k}),
\end{align*}
or 
\begin{equation}\label{PN2}
I_{\cls_k^{\perp}} - \sum_{i=1}^{m} A_{i, k}^* A_{i, k} = (I_{\cls^{\perp}} - \sum_{i=1}^{m} \overline{M}_{z_i}^* \overline{M}_{z_i})_k = X_{k} + \sum_{i=1}^{m} (P_{i, k} - N_{i, k}).
\end{equation}

Here, the subscript $k$ on the middle quantity refers to its restriction to $\cls_k^{\perp}$ which is a reducing subspace. Now,  $$|\mbox{Tr}(I_{\cls_k^{\perp}} - \sum_{i=1}^{m} A_{i,k}^* A_{i,k})| = | \mbox{Tr} (I_{\cls_k^{\perp}} - \sum_{i=1}^{m} A_{i,k} A_{i,k}^*)| = | \mbox{Tr} (I_{\cls^{\perp}} - \sum_{i=1}^{m} \bar{M}_{z_i} \bar{M}_{z_i}^*)_k| = |\mbox{Tr} X_k|.$$
And so, $$|\mbox{Tr}(I_{\cls_k^{\perp}} - \sum_{i=1}^{m} A_{i,k}^* A_{i,k})| = |\mbox{Tr} X_k| \leq \|X_k\| \mbox{dim} \cls_{k-1}^{\perp} = M_0 \|X_k\|.$$
Therefore, $$|\mbox{Tr}(X_{k} + \sum_{i=1}^{m} (P_{i,k} - N_{i,k}))| \leq M_0 \|X_k\|,$$ 
and hence
$$\mbox{Tr}(\sum_{i=1}^{m} P_{i,k}) \leq M_0 (2 \|X_k\| + \sum_{i=1}^{m} \|N_{i,k}\|).$$
But, $\mbox{lim}_{k \raro \infty} (2 \|X_k\| + \sum_{i=1}^{m} \|N_{i,k}\|) = 0$ since $X$ and $N$ are in $\cll^p$, which implies that $\mbox{lim}_{k \raro \infty} \sum_{i=1}^{m} \mbox{Tr} P_{i,k} = 0$. Since $P_{i,k} \geq 0$ we have $\mbox{lim}_{k \raro \infty} \|P_{i,k}\| = 0$. Finally, we have, in view of equation (\ref{PN2}) that $(\overline{M}_{z_1}, \overline{M}_{z_2}, \ldots, \overline{M}_{z_m})$ is an essentially spherical isometry, which concludes the proof. \qed

\vspace{0.2in}
The argument is related to the proof in Section 3 of \cite{GW} for the case of the closure of homogeneous ideals for the $m$-shift space for which the Hilbert-Samuel polynomial is linear; that is, $p_{\cls_k^{\perp}} (k) = M_0 + M_1 k$ for $k \geq K$ and natural numbers $K, M_0$ and $M_1$. Extending the above result to submodules for which the Hilbert-Samuel polynomial has higher degree using this approach would require more control over the $\cll^p$ norms which are not ``linear'' as is the trace.

\vspace{0.6in}

\newsection{Positive regular Hilbert modules}

\vspace{0.4in}

In this section we relate the Hilbert module, $\clh_P$, obtained from a given positive regular polynomial $P$ to the $m$-shift Hilbert module and a certain quasi-homogeneous ideal $J_P$ in $\mathbb{C}[\bm{z}]$. We show that $\clh_P$ is essentially reductive if the corresponding $[J_P]$ is essentially reductive. Here, the number of variables in the $m$-shift space, $H^2_m$, is much greater than that in $P(\bm{z})$. By virtue of Theorem \ref{eqv-question}, it will follow that an affirmative solution showing essential reductivity for homogeneous ideals in Hilbert module, satisfying ($\bm{A}$), ($\bm{B}$) and ($\bm{C}$) will show that all $\clh_P$ are essentially reductive. 

A polynomial 
\begin{equation}\label{P}
P(\bm{z}) = \sum_{i=1}^{m} a_i z_i + \sum_{i=m+1}^{m+K} a_i \bm{z}^{\bm{\alpha}_i},
\end{equation} 
is {\it positive regular} if $a_i >0$ for $i=1, 2, \ldots, m$ and $a_i \geq 0$ for $i = m+1, m+2, \ldots, m+K$, with $\bm{\alpha}_i$ in $\mathbb{N}^m$ (see \cite{BS}, \cite{P} for more on positive regular polynomials). 

We set  
\begin{equation}\label{Pdomain}
\cld_P = \{ \bm{z} \in \mathbb{C}^m : \sum_{i=1}^{m} a_i |z_i|^2 + \sum_{i=m+1}^{m+K} a_i |\bm{z}^{\bm{\alpha}_i}|^2 < 1 \},
\end{equation}

\NI which is a Reinhardt domain in $\mathbb{C}^m$; that is, a domain for which $(z_1, z_2, \ldots, z_m)$ is in $\cld_P$ if and only if $(e^{i \theta_1} z_1, e^{i \theta_2} z_2, \ldots, e^{i \theta_m} z_m)$ is also in $\cld_P$, for $e^{i \theta_j}$ in $\mathbb{T}$, for $j=1, 2, \ldots, m$.

We define the kernel function $k_P$ on $\cld_P \times \cld_P$ so that 
\begin{equation}
k_P(\bm{z}, \bm{w}) = (1 - \sum_{i=1}^{m} a_i z_i \bar{w}_i + \sum_{i=m+1}^{m+K} a_i \bm{z}^{\bm{\alpha}_i} \bar{\bm{w}}^{\bm{\alpha}_i})^{-1},
\end{equation}
for $\bm{z}$ and $\bm{w}$ in $\cld_P$. Let $\clh_P$ be the corresponding reproducing kernel Hilbert space for $k_P$. Then $\clh_P$ is a Hilbert module over $\mathbb{C}[\bm{z}]$, with $\bm{z} = (z_1, z_2, \ldots, z_m)$. If $\delta_{\bm{\alpha}}$ is the coefficient for $\bm{z}^{\bm{\alpha}}$ in the Taylor series  expansion of the function $$\bm{z} \raro (1 - \sum_{i=1}^{m} a_i z_i  + \sum_{i=m+1}^{m+K} a_i \bm{z}^{\bm{\alpha}_i})^{-1},$$ then one can show that 

\begin{equation}
\clh_P = \{ f \in \clo(\cld_P) : f(\bm{z}) = \sum_{\bm{\beta} \in \mathbb{N}^m} b_{\bm{\beta}} \bm{z}^{\bm{\beta}},
\; \mbox{with} \; \sum_{\bm{\beta} \in \mathbb{N}^m} \frac{|b_{\bm{\beta}}|^2}{\delta_{\bm{\beta}}} < \infty \}. 
\end{equation}

Note that, $\|\bm{z}^{\bm{\beta}}\|^2 = \frac{1}{\delta_{\bm{\beta}}}$ and, the $\{\bm{z}^{\bm{\beta}}\}$ are orthogonal and hence  $\clh_P$ satisfies ($\bm{A}$). 

An $m$-tuple of bounded operators $(T_1, T_2, \ldots, T_m)$ on $\clh$ is said to be $P$-{\it contractive} if 
\begin{equation}
 \sum_{i=1}^{m} a_i T_i T_i^* + \sum_{i=m+1}^{m+K} a_i T^{\bm{\alpha}_i} T^{*{\bm{\alpha}_i}} \leq I_{\clh},
\end{equation}

\NI where $T^{\bm{\alpha}} = T_1^{\alpha_1} T_2^{\alpha_2} \cdots T_m^{\alpha_m}$.  
In case $P(\bm{z}) = \sum_{i=1}^{m} z_i$, then a $P$-contractive commuting $m$-tuple $(T_1, T_2, \ldots, T_m)$ is what is often called a row contraction or $\sum_{i=1}^{m} T_i T_i^* \leq I$. Note that $(T_1, T_2, \ldots, T_m)$ is a row contraction if and only if $(T_1^*, T_2^*, \ldots, T_m^*)$ is a spherical contraction. We use both notions in this paper to conform to the literature. 

The following lemmas are essentially from \cite{BS}. 
\vspace{0.2in}

\begin{Lemma}\label{Pcontraction}
For $P$ a positive regular polynomial, one has 
\begin{equation}\label{Pequation}
P_0 = I_{\clh_P} - \sum_{i=1}^{m} a_i M_{z_i} M_{z_i}^* + \sum_{i=m+1}^{m+K} a_i M_{\bm{z}}^{\bm{\alpha}_i} M_{\bm{z}}^{*{\bm{\alpha}_i}},
\end{equation}
where $P_0$ is the orthogonal projection onto the one-dimensional subspace of constant functions and $M_{\bm{z}}^{\bm{\alpha}} = M_{z_1}^{\alpha_1} M_{z_2}^{\alpha_2} \cdots M_{z_m}^{\alpha_m}$.
\end{Lemma}

\vspace{0.2in}

\begin{Lemma}\label{univP}
Let $P$ be a positive regular polynomial and $(T_1, T_2, \ldots, T_m)$ be a commuting $P$-contractive $m$-tuple on a Hilbert space $\clh$ such that there exists a vector $\nu$ in $\clh$ such that $T^{\bm{\beta}} \nu \perp \nu$ for $\beta \neq \bm{0}$. Then the mapping $X_P z^{\bm{\beta}} = T^{\bm{\beta}} \nu$ extends to a contractive module map $X_P : \clh_P \raro \clh$. 
\end{Lemma}

\vspace{0.3in}

Note that if $\alpha_i = 1$ for $i=1, 2, \ldots, m$ and $K=0$, then $\clh_P = H^2_m$, which we know is essentially reductive. We are concerned with the essential reductivity of $\clh_P$ for an arbitrary positive regular polynomial $P(\bm{z})$. 

Let $T = (T_1, T_2, \ldots, T_m)$ be an $m$-tuple of operators on the Hilbert space which is a row contraction; that is, $\sum_{i=1}^{m} T_i T_i^* \leq I_{\clh}$. We define the completely positive map $\eta_T : \cll(\clh) \raro \cll(\clh)$ by $\eta_T (X) = \sum_{i=1}^{m} T_i X T_i^*$, for $X$ in $\cll(\clh)$. Then $$I_{\clh} \geq \eta_T (I_\clh) \geq \eta_T^2 (I_\clh) \geq \cdots,$$ 
and $T$ is called {\it pure} if $A_T = \lim_{n \raro \infty} \eta_T^n (I_\clh) = 0$.

The following theorem is due to Arveson (see \cite{A2}).

\vspace{0.1in}

\begin{Theorem}\label{model}
If $(T_1, T_2, \ldots, T_m)$ is a pure commuting row contractive $m$-tuple acting on a Hilbert space $\clh$, then there exists a Hilbert space $\cle$ and a submodule $\cls$ of $H^2_m \otimes \cle$ so that $\clh \cong H^2_m \otimes \cle/\cls$ as Hilbert modules over $\mathbb{C}[\bm{z}]$ with $T_i \raro \overline{M}_{z_i}$ for $i=1, 2, \ldots, m$, where $\overline{M}_{z_i}$ denotes the image of $M_{z_i} \otimes I_{\cle}$ in $H^2_m \otimes \cle/\cls$. Moreover, $\cle$ can be identified with $\overline{\mbox{Ran}} (I_{\clh} - \sum_{i=1}^{m} T_i T_i^*)$.  
\end{Theorem}
\vspace{0.2in}

For the positive regular polynomial $P(\bm{z})$, consider the $(m+K)$-tuple $$M_P = (\sqrt{a_1} M_{z_1}, \sqrt{a_2} M_{z_2}, \ldots, \sqrt{a_m} M_{z_m}, \sqrt{a_{m+1}} M^{\bm{\alpha}_{m+1}}_{\bm{z}}, \ldots, \sqrt{a_{m+K}} M_{\bm{z}}^{\bm{\alpha}_{m+K}})$$
on $\clh_P$. Since the $m$-tuple $(M_{z_1}, M_{z_2}, \ldots, M_{z_m})$ is $P$-contraction on $\clh_P$, it follows that the  $(m+K)$-tuple $M_P$ is a row contraction. Moreover, since $M_{\bm{z}}^{\bm{\alpha}} 1 \perp 1$ for $\bm{\alpha} (\neq \bm{0})$ in $\mathbb{N}^m$, it follows from Theorem \ref{model} that $\clh_P \cong H^2_{m+K}/\cls_P$ for some submodule $\cls_P$ of $H^2_{m+K}$, since in this case, by equation \ref{Pequation}, the defect $\cle$ is one-dimensional. 

Let us be more precise. Let $Z_1, Z_2, \ldots, Z_{m+K}$ be the variables in $H^2_{m+K}$. Consider the operator $X_P : H^2_{m+K} \raro \clh_P$ defined by $X_P Z_i = \sqrt{a_i} z_i$ for $i=1, 2, \ldots, m$ and $X_P Z_i = \sqrt{a_i} \bm{z}^{\bm{\alpha}_i}$, for $i = m+1, m+2, \ldots, m+K$. Then Lemma \ref{univP} implies that $X_P$ extends to a contractive module map with null space a submodule $\cls_P$ of $H^2_{m+K}$, so that the following diagram 

\setlength{\unitlength}{3mm}
\begin{center}
\begin{picture}(20,14)(0,0)
\put(1,3){$H^2_{m+K}/ \cls_P$} 
\put(7,5.8){$ \bar{X}_P$}
\put(1,10){$ H^2_{m+K}$} \put(10,10){$ \clh_P$}
\put(5.6,11){$ X_P$}
\put(3.4,4.4){ \vector(1,1){5.4}} \put(3.5,10.5){ \vector(1,0){6}}
\put(2.4,9){ \vector(0,-1){4.7}} 

 \end{picture}
 \end{center}

\NI is commutative; that is, the quotient map $\bar{X}_P : H^2_{m+K} \raro H^2_{m+K}/ \cls_P$ is an isometric isomorphism. 

The next question concerns determining the submodule $\cls_P$ concretely. We observe first that $\cls_P$ contains the polynomials $Q_i (\bm{Z}) = Z_i - \lambda_i \bm{Z}^{\bm{\alpha}_i}$ for $i= m+1, m+2, \ldots, m+K$, where $\bm{Z}^{\bm{\alpha}_i} = Z_1^{\alpha_1} Z_2^{\alpha_2} \cdots Z_m^{\alpha_m}$, and $\lambda_i = a_i^{\frac{1}{2}}/ (a_1^{\alpha_1} a_2^{{\alpha_2}} \cdots a_m^{{\alpha_m}})^{\frac{1}{2}}$ for $i = m+1, m+2, \cdots, m+K$. Let $J_P$ denote the ideal in $\mathbb{C}[\bm{Z}]$ generated by $$\{Q_{m+1}(\bm{Z}), Q_{m+2}(\bm{Z}), \ldots, Q_{m+K}(\bm{Z})\}.$$

\vspace{0.3in}

\begin{Theorem}
Given a positive regular polynomial $P(\bm{z})$, we have that $$J_P = \cls_P \cap \mathbb{C}[\bm{Z}].$$ Moreover, $J_P$ is quasi-homogeneous for the weight $\bm{n} = (1, 1, \ldots, 1, |\alpha_{m+1}|, |\alpha_{m+2}|, \ldots, |\alpha_{m+K}|)$ in $\mathbb{N}^{m+K}$ and $\cls_P$ is the closure of $J_P$. 
\end{Theorem}

\NI \textsf{Proof.} We observe that each $Q_i(\bm{Z})$ is quasi-homogeneous with weight $\bm{n}$ since the weighted degrees of $Z_i$ and $\bm{Z}^{\bm{\alpha}_i}$ are both $|\alpha_i|$ for $i= m+1, m+2, \ldots, m+K$. Therefore, $J_P$ is quasi-homogeneous for the weight $\bm{n}$. Now the monomials in both $H^2_{m+K}$ and $\clh_P$ are orthogonal, so one can define the weighted action $\gamma_{\lambda}^{\bm{n}}$ of $\mathbb{T}$ on each of them and these actions intertwine the map $X_P$, or we have

\setlength{\unitlength}{3mm}
 \begin{center}
 \begin{picture}(20,14)(0,0)
 \put(2,3){$ \clh_P$} \put(10,3){$ \clh_P$}
 \put(5.6,2.2){$\gamma_{\lambda}^{\bm{n}}$}
 \put(.7,6.5){$ X_P$} \put(11,6.5){$X_P$}
 \put(2,10){$ H^2_{m+K}$} \put(10,10){$H^2_{m+K}$}
 \put(5.6,11){$\gamma_{\lambda}^{\bm{n}}$}
 \put(3.5,3.5){ \vector(1,0){6}} \put(3.5,10.5){ \vector(1,0){6}}
 \put(2.4,9.2){ \vector(0,-1){5}} \put(10,9.2){ \vector(0,-1){5}}

 \end{picture}
 \end{center}

This implies that $\cls_P$ is quasi-homogeneous in the sense that it is invariant under the action of $\{\gamma_{\lambda}^{\bm{n}} : \lambda \in \mathbb{T}\}$. Thus, the ideal $\tilde{J}_P = \cls_P \cap \mathbb{C}[\bm{Z}]$ is quasi-homogeneous and $\cls_P$ is the closure of $\tilde{J}_P$ in $H^2_{m+K}$. Moreover, we have $J_P \subseteq \tilde{J}_P$ and our goal is to show equality. 

Consider the quotient module $H^2_m/{\tilde{\cls}_P}$, where $\tilde{\cls}_P$ denotes the closure of $\tilde{J}_P$ in $H^2_{m+K}$. We have the quotient maps $\tilde{X}_P : H^2_{m+K} \raro H^2_{m+K}/{\tilde{\cls}_P}$ and $Y_P : H^2_{m+K}/\cls_P \raro H^2_{m+K}/{\tilde{\cls}_P}$ so that the diagram

\setlength{\unitlength}{3mm}
\begin{center}
\begin{picture}(20,14)(0,0)
\put(1,3){$H^2_{m+K}/\cls_P \cong \clh_P$} 
\put(7,5.8){$Y_P$}
\put(1,10){$ H^2_{m+K}$} \put(10,10){$H^2_{m+K}/{\tilde{\cls}_P}$}
\put(5.6,11){$\tilde{X}_P$} \put(.7, 6.7){$ X_P$}
\put(3.4,4.4){ \vector(1,1){5.4}} \put(3.5,10.5){ \vector(1,0){6}}
\put(2.4,9){ \vector(0,-1){4.7}} 

 \end{picture}
 \end{center}
is commutative. However, the construction of $\clh_P$ is a universal one for commuting $P$-contractive $(m+K)$-tuples. Hence, there exists a contractive map $Z_P$ from $\clh_P$ to $H^2_{m+K}/\cls_P$ which commutes with $Y_P$. This implies that $Y_P$ is one-to-one and hence, $\cls_P = \tilde{\cls}_P$, which concludes the proof.  \qed

\vspace{0.2in}

As a corollary we have the reduction of essential reductivity of positive regular Hilbert modules $\clh_P$ to a similar question for homogeneous modules. 

\begin{Theorem}\label{eqv-question}
If all homogeneous ideals in Hilbert module completions of $\mathbb{C}[\bm{z}]$ satisfying ($\bm{A}$), ($\bm{B}$) and ($\bm{C}$) have essentially reductive closures, then every positive regular Hilbert module is essentially reductive. 
\end{Theorem}

\vspace{0.5in}

\newsection{Concluding remarks}

\vspace{0.5in}

As we have indicated, an affirmative answer to Questions \ref{Q1} or \ref{Q2}  would imply an affirmative answer to the conjecture of Arveson; that is, the module over $\mathbb{C}[\bm{z}]$ defined by the closure of a homogeneous module in $H^2_m \otimes \mathbb{C}^r$ is essentially reductive. Let us provide few more details and state this result formally.

\begin{Theorem}
If the answer to Question 1 for an arbitrary $m$-tuple is affirmative, then Arveson's conjecture is valid for $H^2_m \otimes \mathbb{C}^r$; that is, the closure of a homogeneous submodule $S$ in $\mathbb{C} [z_1, z_2, \ldots, z_m] \otimes \mathbb{C}^r$ is essentially reductive for $r$ in $\mathbb{N}$.
\end{Theorem}

\NI \textsf{Proof.} Clearly the restriction of the $m$-tuple $(M_{z_1} \otimes I_{\mathbb{C}^r}, M_{z_2} \otimes I_{\mathbb{C}^r}, \ldots, M_{z_m} \otimes I_{\mathbb{C}^r})$ to $[S]$ is an essentially spherical isometry. Moreover, it follows that $[S]$ is finitely generated using the homogeneity of $S$ and the circle group action. Further, the column operator 
\[ X_{\bm{\lambda}} =  \left( \begin{array}{c}
M_{z_1} \otimes I_{\mathbb{C}^r} - \lambda_1 \\
.\\
.\\
.\\
M_{z_m} \otimes I_{\mathbb{C}^r} - \lambda_m \end{array} \right)\] 

\NI has a finite dimensional kernel and closed range for $\bm{\lambda} = (\lambda_1, \lambda_2, \ldots, \lambda_m)$ in $\mathbb{B}^m$. Therefore, the restriction of $X_{\bm{\lambda}}$ to $[S]$ which yields the operator
\[ \left( \begin{array}{c}
T_{z_1} - \lambda_1 \\
.\\
.\\
.\\
T_{z_m} - \lambda_m \end{array} \right), \] 

\NI has closed range and finite dimensional kernel, where $T_i = (M_{z_i} \otimes I_{\mathbb{C}^r})|_{[S]}$ for $i=1, 2, \ldots, m$. Therefore, the $m$-tuple $(T_1 - \lambda_1, T_2 - \lambda_2, \ldots, T_m - \lambda_m)$ is left semi-Fredholm for $\bm{\lambda}$ in $\mathbb{B}^m$. Hence, the affirmative answer to Question 1 yields the desired result. \qed

\vspace{0.2in}

A variant of Question 1, which we will call Question $1^*$, replaces the assumption that $(M_{z_1} - \lambda_1, M_{z_2} - \lambda_2, \ldots, M_{z_m} - \lambda_m)$ is left semi-Fredholm for $\bm{\lambda}$ in $\mathbb{B}^m$ by the assumption that the essential Taylor spectrum of $(M_{z_1}, M_{z_2}, \ldots, M_{z_m})$ is contained in $\partial \mathbb{B}^m$. It is easy to see that the second assumption implies the first one. To show that an affirmative answer to this modified question implies the Arveson conjecture, one would need to establish that $\sigma_e^{\mbox{Tay}} (T_1, T_2, \ldots, T_m) \subseteq \partial \mathbb{B}^m$ for $T_i = (M_{z_i} \otimes I_{\mathbb{C}^r}) |_{[S]}$ for $i = 1, 2, \ldots, m$, where $S$ is a homogeneous submodule of $H^2_m \otimes \mathbb{C}^r$. This fact would be implied by the validity of the Arveson conjecture but may be easier to prove and hence could be a ``stepping stone'' to that result. 

There is another interesting question that arises when one studies Question $1^*$ by considering the image of the $(T_1, T_2, \ldots, T_m)$ on $\clh$ in the Calkin algebra. More precisely, intuition based on the one-variable case might suggest that the Taylor spectrum of the restriction of a commuting $m$-tuple of normal operators to a proper invariant subspace must be larger, at least if the joint spectrum for the normal operators is nice. However, a family of examples in \cite{EP}, where the joint spectrum is contained in $\partial \mathbb{B}^m$ for $m \geq 3$, shows that this hope is false. (The examples are based on an earlier set of examples due to Izzo \cite{I} showing the failure of polynomial approximation on polynomially convex subsets of the unit sphere in $\mathbb{C}^m$ for $m \geq 3$. Some further results by Izzo, prompted by this topic, appear in \cite{I2}.)

For our application, the normal operators can be assumed to be circularly symmetric which is equivalent to ($\bm{A}$) holding for the closure of the polynomials in the $L^2$ space. We need the result under the assumption that implies ($\bm{A^*}$). Rather than being that precise, for our purpose, however, we focus on circular symmetry.

\begin{Question}
Let $N_1, N_2, \ldots, N_m$ be an $m$-tuple of commuting circularly symmetric normal operators for which $N_1^* N_1 + N_2^* N_2 + \cdots + N_m^* N_m = I_{\clh}$ or, equivalently, $\sigma^{\mbox{Tay}} (N_1, N_2, \ldots, N_m) \subseteq \partial \mathbb{B}^m$. Let $\clm$ be a subspace of $\clh$ invariant for the $N_i$ for which $\clh$ is the minimal reducing subspace containing $\clm$. If $\sigma^{\mbox{Tay}} (T_1, T_2, \ldots, T_m) \subseteq \partial \mathbb{B}^m$, where $T_i = N_i|_{\clm}$ for $i=1, 2, \ldots, m$, does it follow that $\clm = \clh$?
\end{Question}
\vspace{0.2in}

An affirmative answer to the question with the ($\bm{A^*}$) assumption would show that Question $1^*$ has an affirmative answer.

It is not clear whether the coordinate multipliers in the examples in \cite{GRS} are essentially normal. An affirmative answer to Question \ref{Q1} would imply that but the lack of essential normality would provide a counter example. 

To consider the more refined results involving $\cll^p$ commutators, one would  need to provide an affirmative answer to the analogous questions involving the notion of a $p$-essential isometry. 

In connection with Question \ref{Q2}, if one assumes that the Hilbert module is spherically isometric, then one can show that the trace of the restriction of $\sum_{i=1}^{m} [T_i^*, T_i]$ to $\clh_k$ satisfies
\begin{equation}
\mbox{Tr} \{ (\sum_{i=1}^{m} [M_{z_i}^*, M_{z_i}])|_{\clh_k} \} = \binom{m+k-1}{k-1} - \binom{m+k-2}{k-2},
\end{equation}

\NI for $k \geq 1$. As a result, one has $\sum_{i=1}^{m} [M_{z_i}^*, M_{z_i}]$ in $\cll^p$ implies that $p>m$. Unfortunately, these calculations are too crude to establish any positive implication. In particular, if one consider the example in Section 1 involving the Hardy space over the polydisk, one sees that it is impossible to establish that the sum of the self commutators is in $\cll^p$ without making some further assumption such as the left semi-Fredholmness assumption introduced there. On the other hand, if we assume that the eigenvalues of $\sum_{i=1}^{m} [M_{z_i}^*, M_{z_i}]|_{\clh_k}$ are all equal (which is not possible), then the sum of the commutators would be in $\cll^p$ for $p>m$. Thus, resolving the issue is a matter of understanding better the distribution of these eigenvalues. 

An interesting question concerns the converse to the implication of an affirmative answer to Question 3 to the Arveson conjecture. In particular, if one knows the result for the closure of homogeneous ideals in $H^2_m$, does that imply a positive answer for Question \ref{Q3} in which one assumes in addition that the $m$-tuple is a spherical isometry or a spherical contraction?

For such a spherical isometry $(T_1, T_2, \ldots, T_m)$ we have $(T_1^*, T_2^*, \ldots, T_m^*)$ is a spherical contraction. Moreover, if $\clh$ has a wandering cyclic vector, then there is a homogeneous submodule $\cls$ of $H^2_m$ so that $\clh \cong H^2_m/\cls$ by \cite{A2}. Hence, Question \ref{Q3} has an affirmative answer in this case if Arveson's conjecture is valid. The same would be true for the $\cll^p$-analogues. 

In \cite{D}, the first author refined the conjecture of Arveson to state that the quotient module is $p$-essentially reductive for $p > \mbox{dim} \clz$ or the degree of the Hilbert-Samuel polynomial. The results in this paper have no implication for the study of that question.

\end{document}